\documentclass{amsart}
\usepackage{amsmath,amssymb,amsthm}
\theoremstyle{definition}

\theoremstyle{remark}

\numberwithin{equation}{section}

\begin{document}

\title{One more counterexample on sign patterns}

\author{Yaroslav Shitov}

\email{yaroslav-shitov@yandex.ru}


\begin{abstract}
The \textit{sepr-sequence} of an $n\times n$ real matrix $A$ is $(s_1,\ldots,s_n)$, where $s_k$ is the subset of those signs of $+,-,0$ that appear in the values of the $k\times k$ principal minors of $A$. The $12\times 12$ matrix
$$\left(\begin{array}{cccccc|ccc|ccc}
0&0&0&0&0&0&0&0&0&a_1&0&0\\
0&0&0&0&0&0&0&0&0&0&a_2&0\\
0&0&0&0&0&0&0&0&0&0&0&a_3\\
0&0&0&0&0&0&0&0&0&0&0&a_4\\
0&0&0&0&0&0&0&0&0&0&0&a_5\\
0&0&0&0&0&0&0&0&0&0&0&a_6\\\hline
b_1&b_2&0&0&0&0&0&0&0&0&0&0\\
b_3&b_4&0&0&b_5&-b_6&0&0&0&0&0&0\\
0&b_7&b_8&-b_9&b_{10}&b_{11}&0&0&0&0&0&0\\\hline
0&0&0&0&0&0&c_1&0&0&0&0&0\\
0&0&0&0&0&0&0&c_2&0&0&0&0\\
0&0&0&0&0&0&0&0&c_3&0&0&0
\end{array}\right)$$
does always have $s_k=\{0,+,-\}$ if $k=3,6,9$ and $s_k=\{0\}$ otherwise, provided that the variables are positive. However, every principal $9\times 9$ minor that is not identically zero can take values of both signs.
\end{abstract}

\maketitle


\thispagestyle{empty}

This is a counterexample to Conjecture 3.1 in~\cite{HLPD}. The claims can be checked with what I hope is a user friendly \textit{Wolfram Mathematica} worksheet~\cite{myseprfile}.

\end{document}